\documentclass[11pt]{article}

\input{epsf}

\usepackage{amsmath,amssymb} 

\newcommand{\qed}{\hbox{\rule{6pt}{6pt}}}

\newcommand{\Z}{\mathbb{Z}}
\newcommand{\R}{ \mathbb{R} }

\newtheorem{theorem}{Theorem}[section]
\newtheorem{corollary}[theorem]{Corollary}
\newtheorem{lemma}[theorem]{Lemma}
\newtheorem{proposition}[theorem]{Proposition}

\newtheorem{definition}[theorem]{Definition}

\newtheorem{example}[theorem]{Example}

\newtheorem{remark}[theorem]{Remark}

\begin{document}

\title{ Homology Theory for the Set-Theoretic
Yang-Baxter Equation and 
 Knot Invariants from  Generalizations of Quandles }

\author{
J. Scott Carter\footnote{Supported in part by NSF Grant DMS 
\#0301095.}
\\University of South Alabama \\
Mobile, AL 36688 \\ carter@jaguar1.usouthal.edu
\and 
Mohamed Elhamdadi
\\ University of South Florida
\\ Tampa, FL 33620  \\ emohamed@math.usf.edu
\and 
Masahico Saito\footnote{Supported in part by NSF Grant DMS 
\#0301089.}
\\ University of South Florida
\\ Tampa, FL 33620  \\ saito@math.usf.edu
}

\maketitle

\vspace{10mm}

\begin{abstract}
A homology theory is developed for set-theoretic Yang-Baxter 
equations, and knot invariants are constructed by generalized 
colorings by biquandles and Yang-Baxter cocycles.
\end{abstract}

\vspace{10mm}

\section{Introduction}

The introduction of virtual knots by Kauffman~\cite{KauA,KauB}
re-focused attention on algebraic structures that are defined via
diagrams. 
Advantages of using virtual knots were observed for 
the bracket polynomial \cite{KauB} and Vassiliev invariants \cite{GPV}.
The fundamental (Wirtinger)
 groups of virtual knots were studied \cite{Kim,SilWil}
and their geometric interpretations were given \cite{KK}.
Generalizations of Alexander polynomials were studied in relations to
virtual knots \cite{Saw,DanSusan}.
The theory of racks 
and quandles 
(in particular the homology
theory thereof) as exposed in 
\cite{FR,FRS1,FRS2} was used to define  state-sum invariants
for knotted surfaces, as well as for classical and virtual knots~\cite{CJKLS}. 
A generalization of quandles,
called biquandles,  is proposed in \cite{KR}. 
Examples include a generalized Burau matrix used in
\cite{Saw} and \cite{DanSusan}.
The set theoretic solutions to the Yang-Baxter equations are
 studied in detail in the papers \cite{ESG,ESS,LYZ,Solo}. 
Their affine solutions appear among our birack matrices. 
Some of these solutions also appeared in \cite{DanSusan} and \cite{Diamond}.

In this paper a
 homology theory for the YBE is constructed, and cocycles are used
to define knot invariants via colorings of 
(virtual) 
knot diagrams by biquandles and
a  state-sum  formulation.

This paper is organized into two sections; 
Section~\ref{homsec}  develops 
algebraic theories,  and 
Section~\ref{invsec} gives 
applications to knot theory.
In Section~\ref{homsec}, the Yang-Baxter sets are reviewed 
and colorings of cubical complex by Yang-Baxter sets are considered
in Section~\ref{colorsubsec}, for developing homology theories 
for the Yang-Baxter sets in Section~\ref{homsubsec}.
Obstruction and extension cocycles are investigated in 
Section~\ref{obexsubsec} along the lines of group and 
quandle cohomology theories~\cite{CJKLS}. 
In Section~\ref{invsec}, after a review on biquandles 
and their colorings in Section~\ref{prelimsubsec},
the cocycle invariants are defined using the Yang-Baxter 
cohomology theory, and computations and applications are
given in Section~\ref{compsubsec}.

\section{Homology theory for set theoretic Yang-Baxter equations} 
\label{homsec}

\subsection{Yang-Baxter sets and their colorings} \label{colorsubsec}

Let $R=(R_1, R_2)$ be a solution to the {\it set theoretic Yang-Baxter
equation} on a set $X$, that is, $R: X \times X \rightarrow  X \times X$,
written componentwise as $R(x_1, x_2)=(R_1(x_1, x_2), R_2(x_1, x_2))$
for $x_1, x_2 \in X$, satisfies the relation 
$$(R \times 1)(1 \times R)(R \times 1)=(1 \times R)(R \times 1)(1 \times R) $$
where $1$ denotes the identity map. 
The set $X$ with a solution $R$ to the Yang-Baxter equation,
$(X,R)$, is called a {\it Yang-Baxter set}.
Set-theoretic Yang-Baxter 
equations have been studied in 
 \cite{ESG,ESS,LYZ,Solo}, for example.

\begin{example} \label{YBEexamples} {\rm 
Let $k$ be a commutative ring with $1$ and with units $s$ and $t$, 
such that $(1-s)(1-t)=0$. 
Then   $\displaystyle R= \left[ \begin{array}{cc} 1-s & s \\ t & 1-t \end{array} \right] $ satisfies the YBE. 
In particular, for any commutative ring $F$ with $1$,
let $k=F[s^{\pm 1}, t^{\pm 1}] / (1-s)(1-t)$. 
Then the above defined $R$ satisfies 
the conditions.
More specifically, 
$\Z_q$ becomes a YB set for 
integers  $s$ and $t$, where  $q=|(1-s)(1-t)|$
if  $s$ and $t$  are units in  $\Z_q$. 

Another similar example is constructed from matrices. 
Let $R: X \times X \rightarrow  X \times X $ be defined by 
a matrix 
$\displaystyle R= \left[ \begin{array}{cc} E - Y & Y \\ Z & E-Z \end{array} \right] $
with 
 invertible 
matrices $Y$ and $Z$ such that $YZ=ZY$ and $(E-Y)(E-Z)=0$.
Then  $R$ satisfies the YBE. 
Specifically, 
for any commutative 
ring $k$ with $1$, 
$\displaystyle Y= \left[ \begin{array}{cc} 1 & s \\ 0 & 1 \end{array} \right] $
and 
$\displaystyle Z= \left[ \begin{array}{cc} 1 & t \\ 0 & 1 \end{array} \right] $
for any $s, t \in k$ give a solution to the YBE. 
} \end{example}

\bigskip 

Let ${\cal I}_n$ 
 be the $n$-dimensional cube $I^n$ 
regarded as a 
CW (cubical) complex, where $I=[0,1]$ 
and $n$ is a positive integer.
Denote the $k$-skeleton  of ${\cal I}_n$ by ${\cal I}_n^{(k)}$. 
Every $k$-dimensional face of ${\cal I}_n$, for any positive integer $k$,
is another $k$-dimensional cube. 
Give the orientation 
for each $k$-face to be the one 
defined from the order of the coordinate 
axes.   Every $k$-face $\sigma$ is 
regarded as having this orientation fixed, and the same
$k$-face with the opposite orientation is denoted by $-\sigma$. 
In particular, note that every $2$-face can be written as 
$$  \{ \epsilon_1 \} \times \cdots \times  \{ \epsilon_{i-1} \} 
 \times I_i  \times \{ \epsilon_{i+1} \}  \times \cdots \times
  \{ \epsilon_{j-1} \}  \times I_j  \times \{ \epsilon_{j+1} \}
 \times \cdots \times  \{ \epsilon_{n}\}  $$
for some $i$ and $j$, $1 \leq i<  j \leq n$,
where $\epsilon_k = 0$ or $1$, 
and $I_i$, $I_j$ denote the $i$-th,  $j$-th factors of a copy of $I$,
respectively. 
We abbreviate parentheses for simplicity.
Similarly, when the number $0$ or $1$ is placed at the $i$-th factor, 
we denote it by $0_i$ or $1_i$, respectively.

\begin{definition} {\rm
The {\it $($Yang-Baxter$)$ coloring of ${\cal I}_n$ by 
a Yang-Baxter set $(X, R)$} is 
a map $L: E({\cal I}_n) \rightarrow X$, where 
$E({\cal I}_n) $ denotes the set of edges ($1$-face) of ${\cal I}_n$,
with each edge oriented as above, such that 
if 
\begin{eqnarray*}
L(  \epsilon_1  \times \cdots \times  I_i  \times \cdots 
 \times  0_j  \times  \cdots \times   \epsilon_{n}  ) & = & x 
\quad \mbox{and} \\
L(  \epsilon_1 \times \cdots \times   1_i   \times \cdots 
 \times  I_j \times  \cdots \times   \epsilon_{n}  ) & = &y , 
\end{eqnarray*}
then 
\begin{eqnarray*}L(  \epsilon_1  \times \cdots \times   0_i   \times \cdots 
 \times  I_j  \times  \cdots \times   \epsilon_{n}  ) &=&R_1(x,y) \quad
\mbox{and} \\
 L(  \epsilon_1  \times \cdots \times  I_i  \times \cdots 
 \times  1_j  \times  \cdots \times   \epsilon_{n}  ) &=&R_2(x,y) . 
\end{eqnarray*} 
} \end{definition}

\bigskip

\begin{figure}
\begin{center}
\mbox{
\epsfxsize=2in
\epsfbox{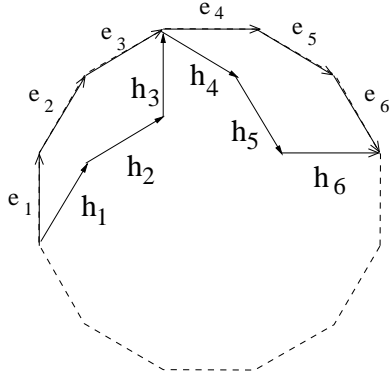} 
}
\end{center}
\caption{ The initial path of ${\cal I}_6$ }
\label{initialpath} 
\end{figure}

To see how many Yang-Baxter colorings  ${\cal I}_n$ admits, 
we specify the {\it initial path} in  ${\cal I}_n$: 
the sequence of edges of  ${\cal I}_n$,
$( e_1, \cdots, e_{n})$, where 
\begin{eqnarray*}
e_1 & = & I_1 \times  0_2 \times \cdots \times  0_n \\
e_2 &=&  1_1  \times I_2 \times  0_3   \times \cdots \times  0_n \\
 \vdots &  \vdots & \quad  \quad  \quad  \quad  \quad  \vdots \\
e_n  &=&   1_1  \times \cdots \times  1_{n-1}   \times I_n ,
\end{eqnarray*}
is called the initial path of   ${\cal I}_n$.
Note that the orientations of the edges $e_i$ 
of the initial path are consistent in the sense that the terminal point
of $e_i$ is the initial point of $e_{i+1}$, for all $i=1, \cdots, n-1$.

For a sub-cube 
$$C=   \epsilon_1   \times \cdots \times 
I_{j_1}  \times \cdots \times 
I_{j_2}  \times \cdots \times 
I_{j_k}  \times \cdots \times   \epsilon_n  $$
of dimension $k$, 
identify $C$ with   ${\cal I}_k$ by 
 the  obvious map
sending the $h$-th factor of $I$ in $I^k$  to $I_{j_h}$ in $C$.
The orientations are preserved by this identification. 
Then the sequence of edges of $C$ corresponding to the initial path
of  ${\cal I}_k$ under this identification is called 
the {\it initial path} of $C$. 
In Fig.~\ref{initialpath}, a 
projection of   ${\cal I}_6$
is depicted, and its initial path are labeled by $e_1, \ldots, e_6$ 
along 
the top edges. 

\begin{lemma} \label{yblemma}
Let $(X,R)$ be a Yang-Baxter set, and $(e_1, \cdots, e_n)$
be the initial path of ${\cal I}_n$.
 For any $n$-tuple of elements of $X$,
$(x_1, \cdots, x_n)$, there exists a unique 
Yang-Baxter coloring $L$ of   ${\cal I}_n$ by $(X,R)$ 
such that $L(e_i)=x_i$ for all $i=1, \cdots, n$. 
\end{lemma}
{\it Proof.\/} 
A sequence of edges 
$(h_1, \ldots, h_n)$ in 
${\cal I}_n$ 
is called  a {\it complete path} 
in ${\cal I}_n$, if 
 the initial point (with respect to the given orientation of each edge)
 of $h_1$ is $0_1 \times \cdots \times 0_n$, 
 the terminal point of $h_n$ is $1_1 \times \cdots \times 1_n$, and 
the initial point 
of $h_{i+1}$ matches the terminal point of $h_{i}$ for all $1<i<n-1$.
For example, the initial path $(e_1, \ldots, e_n)$ is a complete
path, as is the sequence of edges 
$$h_i=e_i, \quad {\mbox{\rm for}} \ \ i\ne j, j+1$$
$$h_j= (\underbrace{1,\ldots, 1}_{j-1},0, I_{j+1}, 
\underbrace{0, \ldots 0}_{n-j-1})$$
$$h_{j+1} = (\underbrace{1,\ldots, 1}_{j-2},0, I_{j},
 \underbrace{0, \ldots 0}_{n-j}) . $$
In Fig.~\ref{initialpath}, (the projection of) another example of a complete
path $(h_1, \ldots, h_6)$ is depicted.

The  complete path $(h_1, \ldots, h_n)$ 
defined above 
is obtained 
from $(e_1, \ldots, e_n)$ by {\it sweeping the square}
 $$(\underbrace{1,\ldots, 1}_{j-1}, I_j \times I_{j+1}, 
 \underbrace{0, \ldots, 0}_{n-j}).$$ 
We describe sweeping  a square 
in full generality as follows. 
Suppose that 
\begin{eqnarray*} 
h_{j} &=&   \epsilon_1  \times \cdots \times  I_{\alpha}   
\times \cdots  \times  0_{\beta} 
\times  \cdots \times   \epsilon_{n} , \\
h_{j+1}  &=&   \epsilon_1 \times \cdots \times   1_{\alpha}   
 \times \cdots  \times  I_{\beta}  
\times  \cdots \times   \epsilon_{n}    , \\
h'_{j}  &=&  \epsilon_1  \times \cdots \times   0_{\alpha} 
 \times \cdots \times  I_{\beta} 
 \times  \cdots \times   \epsilon_{n} ,  \\
h'_{j+1}  &=&  \epsilon_1  \times \cdots \times  I_{\alpha}  
\times \cdots  \times  1_{\beta} 
\times  \cdots \times   \epsilon_{n} ,
\end{eqnarray*}
and $h_i=h'_i$ for $i\neq j,j+1$. 
Then the complete  path  $(h'_1, \ldots h'_n)$
is said to be obtained from the complete  path $(h_1, \ldots h_n)$
by {\it sweeping a square}.
Note that this is not a symmetric relation. 
This situation is also denoted by 
$(h_1, \ldots h_n) \Rightarrow  (h'_1, \ldots h'_n)$.
The two paths $(h_1, \ldots h_n)$ and $(h'_1, \ldots h'_n) $ bound 
the square 
$$ 
(\epsilon_1  \times \cdots \times   I_{\alpha}    \times \cdots \times 
 I_{\beta} \times \cdots  \times   \epsilon_{n} ). $$
This square  is called the {\it square swept by } the paths.
A square swept by paths looks 
like the square in Fig.~\ref{2dboundary},
where $(h_1, \ldots h_n)$ goes through top two edges and $(h'_1, \ldots h'_n)$
goes through the bottom two edges.

If there are paths 
${\cal P}_k= (h_1^{(k)}, \ldots, h_n^{(k)})$ 
for $k=0, \cdots, m$ such that 
${\cal P}={\cal P}_0 \Rightarrow {\cal P}_1 
 \Rightarrow \cdots 
 \Rightarrow {\cal P}_m= {\cal P}'$,
then $ {\cal P}'$ is said to be obtained from ${\cal P}$ by 
{\it sweeping squares}
and is denoted by  $ {\cal P}\Longrightarrow {\cal P}' $ 
(the double arrow in notation is longer in this case). 
By the definition of a Yang-Baxter coloring,
if  the values $L(h_i)$ are specified for all $i=1, \cdots, n$, 
and $(h_i) \Longrightarrow  (h'_j)$,  then there exist values 
$L(h'_j)$  
such that the Yang-Baxter coloring condition is
satisfied for each square that is swept.

We now give a bijective correspondence between the set of complete paths 
and 
the symmetric group 
$\Sigma_n$. 
For a given  
complete path ${\cal P}=(h_1, \ldots, h_n)$, 
define 
a permutation $\sigma=\sigma({\cal P}) \in \Sigma_n$ by 
$$h_{j} =   \epsilon_1  \times \cdots \times  I_{\sigma(j)}   
\times \cdots  
   \epsilon_{n}.$$
That is, 
$\sigma$ is defined so that 
the $j\/$th segment of ${\cal P}$ is parallel to the 
$\sigma(j)\/$th coordinate axis. 

Conversely, given a permutation
$\sigma\in \Sigma_n$, 
we construct 
a unique complete path 
${\cal P}= {\cal P}(\sigma)=(h_1(\sigma), \ldots, h_n(\sigma))$  
such that $\sigma({\cal P})=\sigma$
as follows: 
 \begin{eqnarray*}
h_1(\sigma) & = & (0, \ldots, I_{\sigma(1)}, \ldots, 0), \\
h_2(\sigma) & =& (0, \ldots , \underbrace{1}_{\sigma(1)} , 
\ldots, I_{\sigma(2)},
 \ldots, 0), \\
 & \vdots & \\
h_k(\sigma) & = & (\epsilon_1(\sigma(k)), \ldots, I_{\sigma(k)}, 
\ldots, \epsilon_n(\sigma(k))).
\end{eqnarray*} 
The values of $\epsilon_j(\sigma(k))$ are either $0$ or $1$. 
The value $1$ is achieved in the
$\sigma(1)$ through $\sigma(k-1)$ coordinates. The remaining values are $0$. 
This completes the construction of a bijection.

Now we show that  any complete path ${\cal P}$ 
 can be obtained from the initial path
${\cal P}_0=(e_1, \ldots , e_n)$ by sweeping squares, 
{\it i.e.} 
${\cal P}_0 \Longrightarrow {\cal P}$. 
We prove this inductively by considering
 the permutations $\sigma({\cal P})$ associated to the paths. 
First we introduce some notation on permutations.
Write  the permutation $\left( \begin{array}{cccc} 
1 & 2 & \ldots &n-1 \\ i_1 & i_2 & \ldots &i_{n-1} 
\end{array} \right)\in \Sigma_{n-1}$ as
$(i_1,i_2, \ldots, i_{n-1}).$  
For two permutations $\sigma_1, \sigma_2 \in \Sigma_n$, 
denote by $\sigma_1 \Rightarrow \sigma_2$
if there is a transposition $(i,j)$ such that $\sigma_2 = (i,j)\sigma_1$.
For example, elements of  $\Sigma_3$ have the following relations. 
$$\begin{array}{ccccc} 
(123) &\Rightarrow& (132) & \Rightarrow&  (312)\\
\Downarrow &        & \Downarrow & & \Downarrow \\
(213) &\Rightarrow& (231) & \Rightarrow & (321)
\end{array}$$
By definition, if ${\cal P}_1 \Rightarrow {\cal P}_2$
for two complete paths, then  
$\sigma( {\cal P}_1) \Rightarrow \sigma({\cal P}_2) $. 
Conversely, if $\sigma_1 \Rightarrow \sigma_2$, 
then  ${\cal P}_1 (\sigma_1) \Rightarrow {\cal P}_2( \sigma_2)$.
Indeed, ${\cal P}_2( \sigma_2)$ is obtained from ${\cal P}_1 (\sigma_1)$
by sweeping the square of the form 
$(\epsilon_1, \ldots , I_i, \ldots, I_j, \ldots, \epsilon_n)$,
where  $\sigma_2 = (i,j)\sigma_1$.
Note that 
 every permutation can  be obtained from the identity 
by a sequence of such transpositions. 
It follows that  any complete path ${\cal P}= (h_1, \ldots, h_n)$
is obtained from the initial path $(e_1,\ldots, e_n)$
by sweeping squares: $(e_1,\ldots, e_n)\Longrightarrow (h_1, \ldots , h_n)$.

Also note 
that each edge in the $n$-cube is an edge in some complete path. 
For 
a given
edge $(\epsilon_1, \ldots, I_j, \ldots , \epsilon_n)$ 
consider those values of $k$ for which $\epsilon_k=1$. 
Then connect the given edge back to $(0,0,\ldots, 0)$ by a sequence 
of edges 
parallel to $I_k$. Connect the given edge forward
to $(1,1,\ldots,1)$ by intervals in the remaining coordinates. 
Since every edge is an edge on some complete path, 
and since every complete path can be obtained from the initial 
path by sweeping squares, we have a 
Yang-Baxter
coloring on each edge.

To prove uniqueness, we use projections of knot diagrams, 
and we need the following set-up. 
Consider the path 
$( f_1, \cdots, f_{n})$, where 
\begin{eqnarray*}
f_1 & = & 0_1 \times   \cdots \times 0_{n-1}  \times  I_n \\
f_2 &=&  0_1  \times \cdots \times I_{n-1} \times  1_n \\
 \vdots &  \vdots & \quad  \quad  \quad  \quad  \quad  \vdots \\
f_n  &=&   I_1  \times 1_2 \times \cdots \times  1_{n}   ,
\end{eqnarray*}
which we call the {\it terminal path} of   ${\cal I}_n$.
In Fig.~\ref{initialpath}, the bottom edges form the terminal path,
and the union of the initial and terminal paths form the boundary of 
the  $n$-gon, which is a projection of ${\cal I}_n$ into the plane. 
The union of squares  swept by a sequence 
$(e_i) \Longrightarrow (f_i)$ from the initial path to 
the terminal path is a union of parallelograms tiling the $2n$-gon,
see Fig.~\ref{4dboundary} for such tilings for $n=4$. 
The union of squares swept by any path $(e_i) \Longrightarrow (h_i)$
for any given path $(h_i)$  
is a subset of such a tiling of the whole $2n$-gon. 

Recall the knot diagrams dual to parallelograms as
depicted  in Fig.~\ref{2dboundary}.
Consider the projections of such knot diagrams.
The union of squares swept by a sequence 
$(e_i) \Longrightarrow (f_i)$ gives rise to a projection of
dual knot diagrams that 
consist of generic straight lines on the plane. 
Hence there is a one-to-one correspondence 
between the set of tilings of an $2n$-gon in the plane by parallelograms,
and the set of generic straight $n$ lines on the plane. 
In this correspondence, the Yang-Baxter relation corresponds 
to a Reidemeister type III move, 
as depicted in Fig.~\ref{3dboundary}. 
This is a rearrangement of tilings of $3$ adjacent parallelograms.
We call this a {\it Yang-Baxter rearrangement}.

It remains to be seen that for the tilings 
$T_1$ and $T_2$ corresponding to  given two sequences 
$S_1: (e_i) \Longrightarrow (f_i)$ and 
$S_2: (e_i) \Longrightarrow (f_i)$ respectively,
there is a sequence of Yang-Baxter rearrangement 
$ S_1 \stackrel{a_1}{\rightarrow} \cdots \stackrel{a_k}{\rightarrow} S_2$.

Let $(\ell_1, \cdots, \ell_n)$ be the set of straight lines on the plane 
$\R^2 \times \{0 \}$ 
corresponding to $T_1$, and let $(\ell'_1, \cdots, \ell'_n)$
be those in  $\R^2 \times \{ 1 \}$ 
corresponding to $T_2$. 
Consider
the planes $(P_1, \cdots, P_n)$ in $ \R^2 \times [0,1 ] \subset \R^3$
such that each plane $P_i$ ($i=1, \cdots, n$) is spanned by 
$\ell_i$ and $\ell'_i$.
By a small homotopy if necessary, it is assumed that 
 $(P_1, \cdots, P_n)$ do not have intersections among 
 more than 
three planes, and that the triple points are located at distinct levels
with respect to the height (the second factor of $\R^2 \times [0,1]$) 
direction. 
Such a triple point corresponds to a Reidemeister type III move,
and to a Yang-Baxter rearrangement in the dual tilings.

Let $t_0=0 <  t_1 <  \cdots <  t_v=1$ be numbers such that 
there is exactly one triple point between 
$\R^2 \times \{ t_j \}$ and $\R^2 \times \{ t_{j+1} \}$.
Then the dual tilings of $P_i \cap \R^2 \times \{ t_i \}$
gives a sequence of Yang-Baxter rearrangements from $(\ell_i)$ to 
$(\ell_i')$, as desired.
 \qed

\subsection{Yang-Baxter homology theory} \label{homsubsec}

Let $(X,R)$ be a Yang-Baxter set.
 Let $C_n^{\rm YB}(X)$ be the free 
abelian group generated by
$n$-tuples $(x_1, \dots, x_n)$ of elements of $X$. 

Consider a Yang-Baxter coloring $L$ of  ${\cal I}_n$ 
with $L(e_i)=x_i$ for all $i=1, \cdots, n$, which exists
uniquely by Lemma~\ref{yblemma}. This $L$ is fixed until  the end of 
the definition of the chain complex.
Consider any $k$-face (subcube)  ${\cal J}$ of ${\cal I}_n$. 
Let $(f_1, \cdots, f_k)$ be the initial path of  ${\cal J}$. 
Then there is a unique $k$-tuple of elements of $X$,
$(y_1, \cdots, y_k)$, such that $L(f_j)=y_j$ for all $j=1, \cdots, k$. 
Denote this situation by $L( {\cal J})=(y_1, \cdots, y_k)$.

Let $\partial_n^C$ denote the $n$-dimensional 
boundary map in  cubical homology theory.
Thus $\partial_n^C ({\cal I}_n) = \sum_{i=1}^{2n} \epsilon_i {\cal J}_i$
where 
${\cal J}_i$
is 
an $(n-1)$-face,
$\epsilon_i=\pm 1$ depending on whether the orientation of $ {\cal J}_i$ matches the induced orientation on $ {\cal J}_i$.
For the induced orientation, we take the convention
that the inward pointing normal to an $(n-1)$ face appears
 last in a sequence of vectors that specifies an orientation,
 and the orientation of the $(n-1)$-face is chosen so that 
 this sequence agrees with the orientation of the $n$-cube. 
In particular, $C=I_1 \times \cdots \times I_{n-1}
\times \{ 0 \}$ has 
a  compatible orientation.

Define a
homomorphism
$\partial_{n}: C_{n}^{\rm YB}(X) \to C_{n-1}^{\rm YB}(X)$ by 
$\partial_n ( (x_1, \cdots x_n)) = \sum_{i=1}^{2n} \epsilon_i L(C_i).$
Since $\partial_{n-1}^C \circ \partial_{n}^C =0 $, 
we have $\partial_{n-1} \circ \partial_{n}=0$, 
and $\{ \partial_{n} \}$ defines a chain complex 
$( C_*^{\rm YB}(X),  \partial_{n})$. 
Define as usual the homology groups, 
homology groups with an abelian group coefficient $A$,
cochain groups,
 cohomology groups with an abelian group coefficient $A$, 
and denote them by 
$H_*^{\rm YB}(X)$, $H_*^{\rm YB}(X;A)$, $C^*_{\rm YB}(X;A)$, 
$H^*_{\rm YB}(X;A)$, respectively.
We also call cycles and cocycles in this homology theory 
{\it Yang-Baxter cycles} and {\it cocycles}, respectively.

\begin{definition} {\rm
The homology and cohomology theories defined above are
called the {\it homology and cohomology theories of set theoretic 
Yang-Baxter equation}.
} \end{definition}

\bigskip

We exhibit explicit formulas of the boundary homomorphisms
for low dimensions.

\begin{figure}
\begin{center}
\mbox{
\epsfxsize=2in
\epsfbox{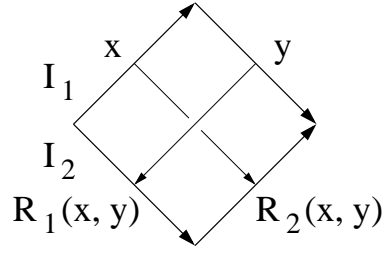} 
}
\end{center}
\caption{ $2$-dimensional boundary homomorphism  }
\label{2dboundary} 
\end{figure}

\begin{example} {\rm 
Let $(X,R)$ be a Yang-Baxter set. 
In Fig.~\ref{2dboundary}, the $2$-dimensional cube ${\cal I}_2$ is depicted. 
The top two edges form the initial path, and the edges are colored 
by $x,y \in X$. The bottom edges are colored 
by $R_1(x,y)$ and $R_2(x,y)$, so that these assignment indeed defines
a Yang-Baxter coloring of a square. 
Thus the boundary homomorphism in this case is given by 
$$ \partial_2 (x,y) = (x) + (y) - (R_1(x,y)) - (R_2(x,y)). $$

In this square, a correspondence between a square and a positive crossing 
point used for  classical knot diagrams is depicted. We use this 
correspondence in the following examples.
} \end{example}

\begin{figure}
\begin{center}
\mbox{
\epsfxsize=4in
\epsfbox{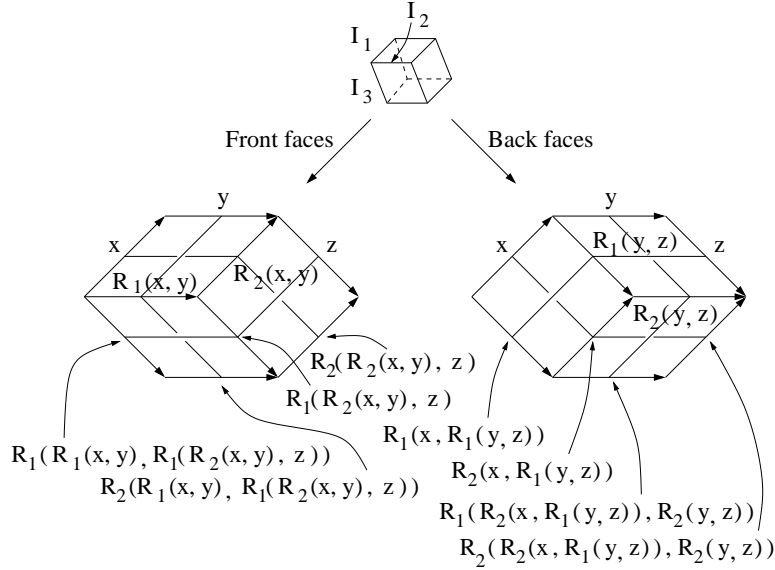} 
}
\end{center}
\caption{ $3$-dimensional boundary homomorphism  }
\label{3dboundary} 
\end{figure}

\begin{example} {\rm 
In Fig.~\ref{3dboundary}, the $3$-dimensional cube ${\cal I}_3$ is depicted
at the top of the figure. 
The top three paths form the initial path, and the elements $x,y,z$ 
are assigned.
The cube, as depicted in the figure, consists of 
three front faces and three back faces from the reader's perspective.
In the top figure of the cube, the three back faces are located 
behind the three front  faces, and thus depicted by dotted lines. 
In the bottom of the figure, the front faces (left) 
and the back faces (right) are depicted separately. 

The front faces determine, via the condition of a Yang-Baxter coloring,
the colors assigned to all the edges of the front faces as depicted. 
The same is true for the edges of the back faces. 
There are six edges that are common for both front and back faces.
Three edges of the initial path (labeled by $x,y,z$) at the top, 
and three edges at the bottom. The expressions of the colors 
assigned to these bottom three edges obtained from front faces and those 
obtained from back faces are different, as depicted in the figure.
These are, of course, the same elements in $X$ for each edge, 
since $R$ is a Yang-Baxter solution, and the equalities that these 
elements are the same indeed are equivalent to the condition that 
$R$ is a solution. The relation between cubes and crossings is again 
depicted in the bottom figure, and it is seen that the Yang-Baxter 
equation corresponds to the type III Reidemeister move, 
as known in knot theory. 
The boundary homomorphism in this case is given by 
\begin{eqnarray*}
\lefteqn{ \partial_3(x,y,z) } \\
&=& (x,y)+(R_2(x,y),z)+ (R_1(x,y), R_1( R_2(x,y),z) ) \\
 &- & \{ (y,z) + (x, R_1(y,z)) + (R_2(x, R_1(y,z)), R_2(y,z)) \} 
\end{eqnarray*} 
In terms of knot diagrams,
it is seen that each colored crossing contributes a $2$-chain (pair)
$(\alpha, \beta)$, and the Reidemeister type III move 
(before $-$ after) gives the boundary map.
} \end{example}

\begin{figure}
\begin{center}
\mbox{
\epsfxsize=4.5in
\epsfbox{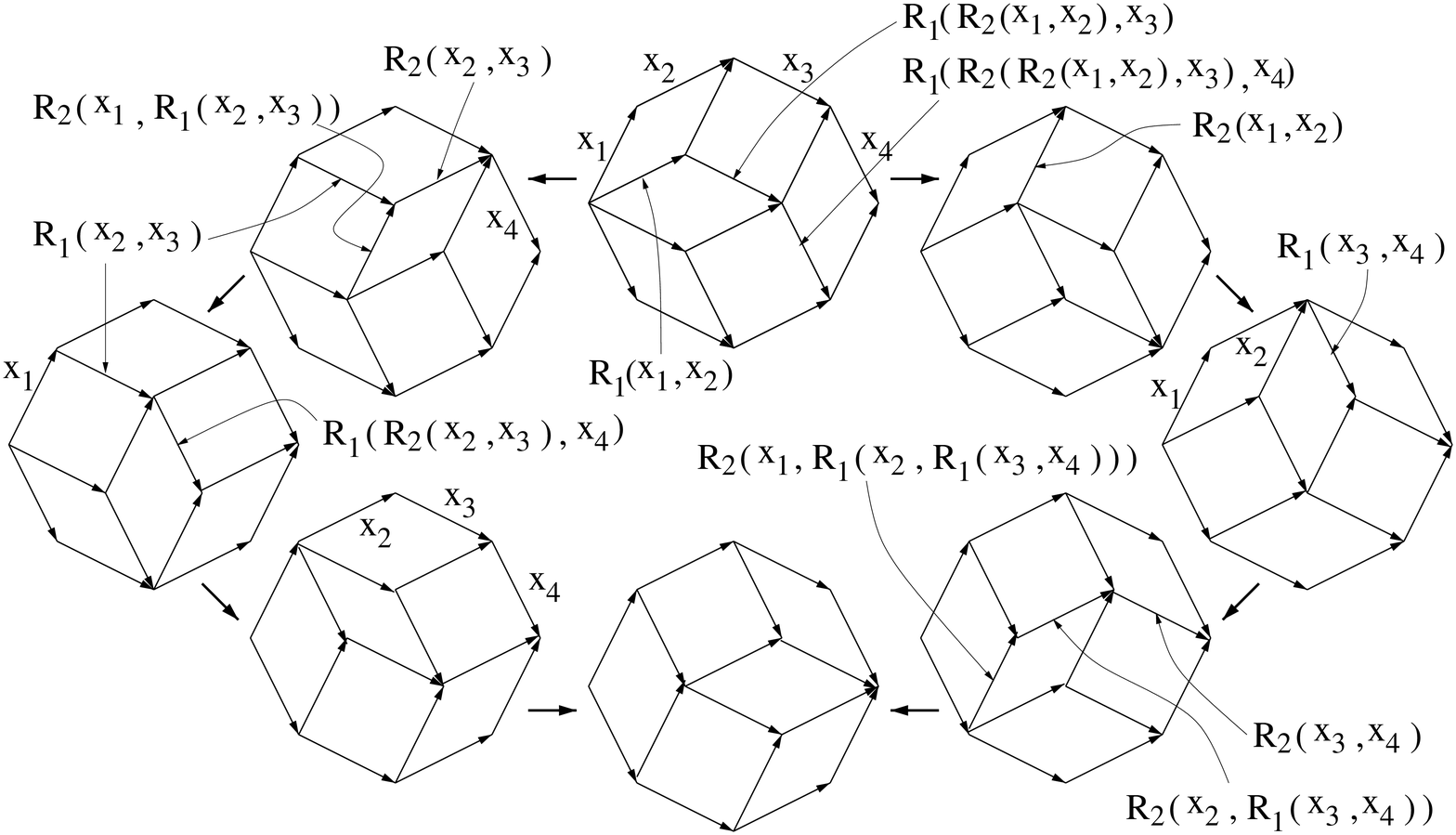} 
}
\end{center}
\caption{ $4$-dimensional boundary homomorphism  }
\label{4dboundary} 
\end{figure}

\begin{example} {\rm 
In Fig.~\ref{4dboundary}, 
the $4$-dimensional case is depicted. 
A  triple 
$(x_1, x_2, x_3) \in C_3^{\rm YB}(X)$ is represented by 
a cube depicted in Fig.~\ref{3dboundary} whose initial path 
$(e_1, e_2, e_3)$ is 
colored by  $(x_1, x_2, x_3)$. As in Fig.~\ref{3dboundary},
such a cube is depicted by the  three front faces (left) and
 three back faces (right). 
In   Fig.~\ref{4dboundary} left and right sides,
the four front and back $3$-faces of ${\cal I}_4$ are depicted respectively,
by exhibiting each $3$-face (cube) by showing the front three and back three.

For example, the change (left to right) of  Fig.~\ref{3dboundary}
happens for a cube with initial path colored by $(x_1, x_2, x_3)$
as the first step in the left of Fig.~\ref{4dboundary} (from top right 
to top middle), contributing the term $(x_1, x_2, x_3)$
in  $\partial_4(x_1, x_2, x_3, x_4) $. 
On the other hand, in the right of Fig.~\ref{4dboundary}, the first change
is from top left to top middle, which is applied to the cube with 
initial path colored by $ (R_1(x_1, x_2), R_1( R_2(x_1, x_2), x_3), 
R_1(R_2( R_2(x_1, x_2), x_3), x_4) )$, the first negative term.
{}From the figure we obtain 
\begin{eqnarray*}
\lefteqn{ \partial_4(x_1, x_2, x_3, x_4) } \\
&=& (x_1, x_2, x_3) + (R_2( x_1, R_1(x_2, x_3) ), R_2( x_2, x_3), x_4) \\
& &+ (x_1,  R_1(x_2, x_3) , R_1(R_2(x_2, x_3) , x_4 ) )+ ( x_2, x_3 , x_4 ) \\
&-& \{ (R_1(x_1, x_2), R_1( R_2(x_1, x_2), x_3), 
R_1(R_2( R_2(x_1, x_2), x_3), x_4) ) 
+ (R_2(x_1, x_2), x_3, x_4 ) \\
& & + ( x_1, x_2, R_1(x_3, x_4)) 
+ (R_2(x_1, R_1(x_2, R_1(x_3, x_4)) ), R_2(x_2, R_1(x_3, x_4)), 
R_2(x_3, x_4) ) \} .
\end{eqnarray*}

\begin{figure}
\begin{center}
\mbox{
\epsfxsize=2.5in
\epsfbox{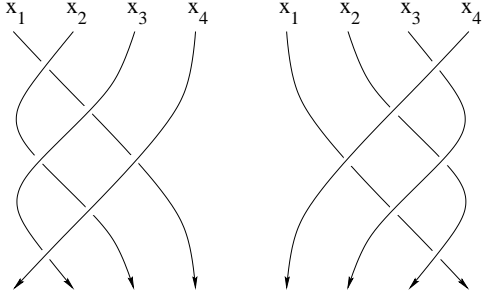} 
}
\end{center}
\caption{Dual knot diagrams }
\label{tetramove} 
\end{figure}

These Fig.~\ref{4dboundary} are dual figures of 
movie version of the ``tetrahedral move,'' one of the Roseman moves
(analogs of Reidemeister moves for knotted surfaces, see \cite{CS:book}).
Using the convention depicted in Fig.~\ref{2dboundary}
on relation between squares and crossings, we obtain the left figure of 
Fig.~\ref{tetramove} from the top  figure of 
Fig.~\ref{4dboundary}. The last figure,
at the bottom  in Fig.~\ref{4dboundary}, 
 corresponds to the right figure in 
Fig.~\ref{tetramove}.
Two changes in  the left and right of 
Fig.~\ref{4dboundary}
correspond to two sequences of type III Reidemeister moves
changing the left figure of Fig.~\ref{tetramove} to the right
(see \cite{CS:book} again). 
In terms of  projections 
into $3$-space 
of a knotted surfaces in $4$-space, 
each type III move corresponds to a triple point in projection.

} \end{example}

\subsection{Obstruction and extension cocycles} \label{obexsubsec}

Low dimensional Yang-Baxter cocycles have interpretations as
obstruction cocycles in extensions of Yang-Baxter sets.
This is a situation similar to group and other cohomology theories,
and generalizes the quandle cohomology case.

Let $0 \rightarrow N \stackrel{i}{\rightarrow}  G \stackrel{p}{\rightarrow} 
A \rightarrow 0 $ 
be an exact sequence of abelian groups 
with a set-theoretic section $s: A \rightarrow G$
which is {\it normalized} in the sense that $s(0)=0$.

Let  $f \in Z^n_{\rm YB}(X;A)$ for a positive integer $n$, 
then $\delta f (x_1, \cdots, x_n)$ has $2n$ terms, each of which 
has $(n-1)$ arguments. Let 
$\delta f (x_1, \cdots, x_n)=T_1 + \cdots + T_{2n}$ be such terms.
Consider the expression $H=s(T_1) + \cdots + s(T_{2n})$.
Then $p(H)=0 \in A$ since $p \circ s = \mbox{id}$, $p$ is a homomorphism,
and $f$ is a cocycle.
Therefore there is a unique  element $\psi(x_1, \cdots, x_n) \in N$ 
such that 
$i \psi(x_1, \cdots, x_n) = H$. 
For example, for $n=2$, we obtain 
\begin{equation}\label{defining2} 
sf(x) + sf(y)= i \psi (x, y) +  sf(R_1(x,y)) + sf(R_2(x,y)).
\end{equation}

\begin{proposition} $\psi \in Z^{n+1}_{\rm YB}(X;N)$. 
\end{proposition}
{\it Proof.\/}
The $2n$ terms $T_1, \cdots, T_{2n}$ are in one-to-one correspondence 
with $(n-1)$-faces of ${\cal I}_n$ whose initial path 
is labeled by $x_1, \cdots, x_n$ via $L$. Assign $i \psi(x_1, \cdots, x_n)$
to this cube ${\cal I}_n$.
Then in ${\cal I}_{n+1}$ whose initial path 
is labeled by $(x_1, \cdots, x_{n+1})$, the $n$-faces ${\cal I}$
are in one-to-one correspondence with the assigned $i\psi$. 
Then $\delta \psi =0$ follows from $\delta^2=0$. 
\qed

\noindent
We call such a cocycle $\psi$ an {\it obstruction 
$(n+1)$-cocycle}.

Explicit calculations can be carried out using this correspondence.
For example, for $n=2$, 
on the one hand,
 \begin{eqnarray*}
\lefteqn{ 
 \underline{sf(x) + sf(y)} + sf(z)} \\ &=& 
 i \psi (x, y) +  sf(R_1(x,y)) + \underline{sf(R_2(x,y)) + sf(z)} \\
 &=&  i \psi (x, y) +  i \psi (R_2(x,y), z)  + \underline{ sf(R_1(x,y))
 + sf( R_1( R_2(x,y), z) )} + sf( R_2( R_2(x,y), z) ) \\
 &=&  i \psi (x, y) +  i \psi (R_2(x,y), z) +  i \psi(R_1(x,y), R_1( R_2(x,y), z) ) \\ & & 
+  sf(R_1( R_1(x,y),  R_1( R_2(x,y), z)))
+ sf( R_2(  R_1(x,y),  R_1( R_2(x,y), z)) ) ) +  sf( R_2( R_2(x,y), z) )
\end{eqnarray*} 
and on the other hand, 
 \begin{eqnarray*}
\lefteqn{ 
 sf(x) +  \underline{sf(y)+ sf(z)} } \\ &=& 
 i \psi (y,z) + \underline{sf(x) + sf(R_1(y,z))} + sf(R_2(y,z)) \\
 &=&  i \psi (y,z) +  i \psi (x, R_1(y,z) ) +
sf( R_1(x, R_1(y,z)) ) +  \underline{sf( R_2( x, R_1(y,z)))+sf( R_2(y,z)) } \\
 &=&  i \psi (y,z) +  i \psi (x, R_1(y,z) ) +
 i \psi ( R_2( x, R_1(y,z)),  R_2(y,z)) \\ & & +
sf( R_1(x, R_1(y,z)) )
+ sf( R_1(  R_2( x, R_1(y,z)),  R_2(y,z)) ) +
sf( R_2(   R_2( x, R_1(y,z)),  R_2(y,z)) )
\end{eqnarray*} 
where the defining relation of $\psi$ was applied to underlined terms.
We recover the formula for the second coboundary homomorphism this way. 
This computation 
is  directly visualized from Fig.~\ref{3dboundary}.

Next we consider extensions.

\begin{proposition} \label{extprop}
Let $V=A \times X$  
and $S: V \times V \rightarrow  V \times V $
be defined by 
$$S( (a_1, x_1), (a_2, x_2) ) = ( (a_2 + \psi_1(x_1, x_2) , R_1(x_1, x_2) ),
 (a_1 + \psi_2(x_1, x_2),  R_2(x_1, x_2) )$$
for any $(a_i, x_i) \in V$. 
If  $(V, S)$ is a Yang-Baxter set, then 
$\psi_1 + \psi_2 \in  Z^2_{\rm YB}(X;A)$.
\end{proposition}
{\it Proof.\/}
We compute
\begin{eqnarray*}
\lefteqn{ (S \times 1 )( 1 \times S)(S \times 1)(
 ( a_1, x_1), (a_2, x_2), (a_3, x_3) )  } \\
&=& ( \quad ( \quad  a_3 + \psi_1( R_2( x_1, x_2) , x_3) + 
 \psi_1(R_1(x_1, x_2), R_1( R_2(x_1, x_2), x_3) ) , \\
 & & \quad  R_1( R_1(x_1, x_2), R_1(R_2( x_1, x_2), x_3) ) \quad  ),\\
& & (  \quad a_2 + \psi_1(x_1, x_2) + \psi_2( R_1(x_1,x_2), 
R_1( R_2( x_1, x_2), x_3) ) , \\
 & &  \quad  R_2( R_1(x_1, x_2), R_1(R_2( x_1, x_2), x_3) \quad  )  \quad ), \\
& & ( \quad a_1 +  \psi_2(x_1, x_2) + \psi_2( R_2( x_1, x_2) , x_3), \\
& & \quad  R_2(R_2(x_1, x_2) , x_3) \quad  )  \quad )
\end{eqnarray*} 
and on the other hand, 
\begin{eqnarray*}
\lefteqn{( 1 \times S)(S \times 1)( 1 \times S)(
 ( a_1, x_1), (a_2, x_2), (a_3, x_3) )  } \\
&=& ( \quad  ( a_3 + \psi_1( x_2, x_3) + \psi_1( x_1, R_1(x_2, x_3) ) , \\
 \quad & & R_1(x_1, R_1(x_2, x_3)) \quad  ), \\
& & ( \quad a_2 +  \psi_2( x_2, x_3) + \psi_1(R_2(x_1, R_1(x_2, x_3)),   R_2(x_2, x_3) ),\\
& &  \quad R_1(R_2(x_1, R_1(x_2, x_3)),  R_2(x_2, x_3))\quad ), \\
& & (  \quad a_1 + \psi_2(x_1, R_1(x_2, x_3) ) 
+   \psi_2(R_2(x_1, R_1(x_2, x_3)),   R_2(x_2, x_3)) , \\
& & \quad R_2(R_2(x_1, R_1(x_2, x_3)),  R_2(x_2, x_3)) \quad) \quad) .
\end{eqnarray*} 
Hence we obtain 
\begin{eqnarray*}
\lefteqn{ \psi_1( R_2( x_1, x_2) , x_3) + 
 \psi_1(R_1(x_1, x_2), R_1( R_2(x_1, x_2), x_3) } \\
&=& \psi_1( x_2, x_3) + \psi_1( x_1, R_1(x_2, x_3) ) \\
\lefteqn{\psi_1(x_1, x_2) + \psi_2( R_1
(x_1,x_2), 
R_1( R_2( x_1, x_2), x_3)} \\
&=&  \psi_2( x_2, x_3) + \psi_1(R_2(x_1, R_1(x_2, x_3)),   R_2(x_2, x_3) ) \\
\lefteqn{ \psi_2(x_1, x_2) + \psi_2( R_2( x_1, x_2) , x_3) } \\
&=&  \psi_2(x_1, R_1(x_2, x_3) ) 
+   \psi_2(R_2(x_1, R_1(x_2, x_3)),   R_2(x_2, x_3))
\end{eqnarray*} 
{}from each factor containing $a_3, a_2, a_1$ respectively,
and by adding these equalities we obtain  the $2$-cocycle 
condition for $\psi_1+ \psi_2$, and the result follows.
\qed

\begin{corollary} \label{extcor}
Suppose $A$ is a ring in which $2$ is invertible. 
Let $V=A \times X$ and $S: V \times V \rightarrow  V \times V $
be defined by 
$$S( (a_1, x_1), (a_2, x_2) ) = ( (a_2 + \psi(x_1, x_2) , R_1(x_1, x_2) ),
 (a_1 + \psi(x_1, x_2),  R_2(x_1, x_2) )$$
for any $(a_i, x_i) \in V$. 
If  $(V, S)$ is a Yang-Baxter set, then 
$\psi \in  Z^2_{\rm YB}(X;A)$.
\end{corollary}
{\it Proof.\/}
In the preceding proposition, set $\psi_1=\psi_2$, 
then we obtain twice the cocycle condition at the end, 
and obtain the result as $2$ is invertible.
\qed

\begin{definition} \label{extdef} {\rm 
The Yang-Baxter set $(V, S)$ defined in  Proposition~\ref{extprop}
is called an {\it extension} of $(X, R)$ by  $(\psi_1,  \psi_2)$. 

} \end{definition}


\begin{example} \label{omegaextex} {\rm 
Let $X=\Omega_q=\Z_q [s^{\pm 1}, t^{\pm 1} ] / (1-s)(1-t)$.
The same matrix 
 $\displaystyle R= \left[ \begin{array}{cc} 1-s & s \\ t & 1-t \end{array} \right]$ 
 used in Example~\ref{YBEexamples} 
 defines a Yang-Baxter set structure on 
 $\Omega_q$.
Let $A=X=\Omega_q$ and 
define 
$\psi_i: X \times X \rightarrow A(=X)$ by 
$\psi_i(x, y)=u_i(y-x)$ for some $u_i \in X$, for $i=1,2$. 

Let $V=A \times X = (\Omega_q)^2 $ and 
consider $S: V \times V \rightarrow V \times V$
defined  as in the Proposition~\ref{extprop}: 
$$S( (a_1, x_1), (a_2, x_2) ) = ( (a_2 + \psi_1(x_1, x_2) , R_1(x_1, x_2) ),
 (a_1 + \psi_2(x_1, x_2),  R_2(x_1, x_2) )$$
for any $(a_i, x_i) \in V$. 
Then $S$ is represented by the matrix 
$\displaystyle  \left[ \begin{array}{cccc}
0 & -u_1 & 1 & u_1 \\
0 & 1-s & 0 & s \\
1 & -u_2 & 0 & u_2 \\
0 & t & 0 & 1-t 
\end{array} \right] $, also written by $S$. 
Then defining
 $\displaystyle Y= \left[ \begin{array}{cc} 1 & u_1 \\ 0 & s \end{array} \right]$
 and
$\displaystyle Z= \left[ \begin{array}{cc} 1 & -u_2 \\ 0 & t \end{array} \right]$,
$S$ can be written 
 as 
 $\displaystyle S= \left[ \begin{array}{cc} E-Y & Y \\ Z & E-Z \end{array} \right]$
where $Y$, $Z$ invertible.
If $u_2(1-s)=0=u_1(1-t)$, then $YZ=ZY$ and 
 $(E-Y)(E-Z)=0$, 
so that from Example~\ref{YBEexamples}  
$(V,S)$ defines a Yang-Baxter set, and is an extension of $(X,R)$.

} \end{example}

\begin{example} \label{extex} {\rm 
Let $(X, R)$ be the Yang-Baxter set in Example~\ref{YBEexamples}.
Specifically, for a ring $k$
in which $2$ is invertible, let $X=k^2$ and 
$\displaystyle R= \left[ \begin{array}{cc} E-Y & Y \\ Z & E-Z \end{array} \right]$,
where $R$ is a block matrix of $2$ by $2$ matrices 
with $\displaystyle Y=\left[ \begin{array}{cc} 1 & s \\ 0 & 1 \end{array} \right]$
and $\displaystyle Z=\left[ \begin{array}{cc} 1 & t \\ 0 & 1 \end{array} \right]$, where $s, t \in k$. 

Let $A=X=k^2$ and $V=A \times X$. 
Let $\psi : X \times X \rightarrow A (=X)$ be defined by 
$$\displaystyle \psi \left( \left[ \begin{array}{cc} x_1 \\ x_2 \end{array} \right] ,
 \left[ \begin{array}{cc} y_1 \\ y_2 \end{array} \right] \right)
=  \left[ \begin{array}{cc} w (y_2 - x_2)  \\ 0 \end{array} \right] , $$
where $w \in k$. 
Define $S : V \times V \rightarrow  V \times V $
as in the Proposition~\ref{extprop}: 
$$S( (\vec{a}, \vec{x}), (\vec{b}, \vec{y}) ) = ( (\vec{b} + \psi(\vec{x}, \vec{y}) , R_1(\vec{x}, \vec{y}) ), (\vec{a }+ \psi(\vec{x}, \vec{y}),  R_2(\vec{x}, \vec{y}) ) )$$ 
for any $(\vec{a}, \vec{x}), (\vec{b},\vec{y})  \in V$.

We now show that 
$(V,S)$ is an extension of $(X,R)$ by  $\psi$, and as 
a consequence of Corollary~\ref{extcor}, we obtain 
$\psi \in  Z^2_{\rm YB}(X;A)$.

The map $S$ above, on the space $(A \times X) \times (A \times X)$, 
can be written as matrix: 
$$ \left[ \begin{array}{cccc} O & -W & E & W \\
O & E-Y & O & Y \\ E & -W &  O & W \\ O & Z & O & E-Z \end{array} \right] , $$
where $O$ denotes the zero matrix,
and  $\displaystyle W =  \left[ \begin{array}{cc} 0 & w \\ 0 & 0 \end{array} \right]$. 
Then the matrix of $S$ is written as  
$\displaystyle \left[ \begin{array}{cc} E - Y' & Y' \\ Z' & E-Z'  \end{array} \right]$,
where $\displaystyle Y'=  \left[ \begin{array}{cc}E & W \\ O & Y  \end{array} \right]$
and  $\displaystyle Z'=  \left[ \begin{array}{cc}E & -W \\ O & Z  \end{array} \right]$.
It is seen that $(E-Y')(E-Z')=O$,
and
$(V,S)$ is a Yang-Baxter set. 

Corollary~\ref{extcor} implies that $\psi \in Z^2_{\rm YB}(X;A)$
provided $2$ is invertible in $k$.
} \end{example}

\begin{proposition} 
Let $k$ be a ring such that $2$ is invertible, set $X=k^2$, 
and $(X,R)$ be the Yang-Baxter set defined 
in Example~\ref{YBEexamples} and Example~\ref{extex},
with $w \neq 0$ and $s=t$: 
$\displaystyle R= \left[ \begin{array}{cc} E-Y & Y \\ Y & E-Y \end{array} \right]$,
 $\displaystyle Y=\left[ \begin{array}{cc} 1 & t \\ 0 & 1 \end{array} \right]=Z$. 
Then $H^2_{\rm YB} (X;A) \neq 0$. 
\end{proposition}
{\it Proof.\/}
Note that 
$\displaystyle \psi \left( \left[ \begin{array}{cc} 1 \\ 0 \end{array} \right] ,
 \left[ \begin{array}{cc} 0 \\ 1 \end{array} \right] \right)
=  \left[ \begin{array}{cc} w   \\ 0 \end{array} \right] $
and 
$\displaystyle \psi \left( \left[ \begin{array}{cc} 0 \\ 1 \end{array} \right] ,
 \left[ \begin{array}{cc} 1 \\ 0 \end{array} \right] \right)
=  \left[ \begin{array}{cc} - w   \\ 0 \end{array} \right] $,
so that they take distinct values, since $w \neq 0$
and $2$ is invertible in $k$. 

On the other hand, suppose $\psi=\delta f$ for some 
$f \in C^1_{\rm YB}(X;A)$. 
One computes that 
\begin{eqnarray*}
 (\delta f)  \left( \left[ \begin{array}{cc} 1 \\ 0 \end{array} \right] ,
 \left[ \begin{array}{cc} 0 \\ 1 \end{array} \right] \right)
&= &f  \left( \left[ \begin{array}{cc}  1 \\ 0 \end{array} \right] \right)
+  f  \left( \left[ \begin{array}{cc} 0 \\ 1 \end{array} \right] \right)
-  f  \left( \left[ \begin{array}{cc} t \\ 1 \end{array} \right] \right)
- f  \left( \left[ \begin{array}{cc}-t+1 
\\ 0  \end{array} \right] \right)
\\ &=&  (\delta f)  \left( \left[ \begin{array}{cc} 0 \\ 1 \end{array} \right] ,
 \left[ \begin{array}{cc} 1 \\ 0 \end{array} \right] \right), 
\end{eqnarray*}
which is a contradiction.
\qed

\bigskip

We consider another family of examples that are similar to 
those considered in \cite{CENS}.

\begin{proposition} \label{dbleomegaprop} 
For positive integers $q, h, k$, let 
$X=\Omega_q^{(h,k)}= 
\Z_q [s^{\pm 1}, t^{\pm 1} ]/ \{ (1-s)^h, (1-t)^k \}$. 
Let $R: X \times X \rightarrow X \times X $ be represented by 
the matrix $\left[ \begin{array}{cc} 1-s &s\\t&1-t \end{array} \right]$.
Then $V=\Omega_q^{(h+1,k+1)}$ is an extension of $\Omega_q^{(h,k)}$,
where $S:  V \times V \rightarrow V \times V $
is represented by the same matrix 
that defines 
$R$,
and $A=\Z_q \times \Z_q$. 
\end{proposition} 
{\it Proof.\/} 
Recall that $\Omega_q= \Z_q [t^{\pm 1}, s^{\pm 1}]/((1-s)(1-t)).$
Represent elements $\alpha,\beta \in V$ by 
$\alpha= \alpha_0 + \sum_{i=1}^h \alpha_i (1-s)^i +  \sum_{j=1}^k \alpha'_j (1-t)^j $
and $\beta= \beta_0 + \sum_{i=1}^h \beta_i (1-s)^i +  \sum_{j=1}^k \beta'_j (1-t)^j $.
Let $\bar{\alpha}=p(\alpha)$ and $\bar{\beta}=p(\beta)$ where
$p: V \rightarrow X$ is the natural projection. 
Let $\sigma: X \rightarrow V$ be the natural set-theoretic section defined by  
$\sigma( \alpha_0 + \sum_{i=1}^{h-1} \alpha_i (1-s)^i +  \sum_{j=1}^{k-1} \alpha'_j (1-t)^j )
=   \alpha_0 + \sum_{i=1}^h \alpha_i (1-s)^i +  \sum_{j=1}^k \alpha'_j (1-t)^j $
where $\alpha_h = \alpha'_k=0$.
We identify $V$ with $A \times X$ by $f:V \rightarrow A \times X$, where  
$$ f( \alpha_0 + \sum_{i=1}^h \alpha_i (1-s)^i +  \sum_{j=1}^k \alpha'_j (1-t)^j )
= (  (\alpha_h, \alpha'_k) ,  \bar{\alpha}) .$$
Then we compute that 
\begin{eqnarray*} 
\lefteqn{S(\quad ( (\alpha_h, \alpha'_k),\quad \bar{\alpha} ) , \quad ( (\beta_h, \beta'_k), \quad \bar{\beta} )\quad) }\\
&=& S( \alpha, \beta) \\
&=& (\quad (1-s)\alpha + ( 1-(1-s))\beta,\quad (1-(1-t))\alpha + (1-t)\beta \quad) \\
&=& (\quad \beta + (1-s)(\alpha-\beta),\quad \alpha + (1-t) (\beta-\alpha)\quad ) \\
&=& (\quad (\beta_h+ (\alpha_{h-1}-\beta_{h-1})  (1-s)^h + \beta'_k (1-t)^k +
\sigma( S_1( \bar{\alpha} , \bar{\beta}) ) , \\
 & &  \alpha_h (1-s)^h + (\alpha'_k+ (\beta'_{k-1}-\alpha'_{k-1})  (1-t)^k +
\sigma( S_2( \bar{\alpha} , \bar{\beta}) )\quad ) \\
&=& ( \quad ( \quad \beta_h+ (\alpha_{h-1}-\beta_{h-1}), \quad \beta'_k,\quad S_1( \bar{\alpha} , \bar{\beta}) \quad ) , \\ & & 
 ( \quad  \alpha_h,\quad \alpha'_k+ (\beta'_{k-1}-\alpha'_{k-1}),\quad  S_2( \bar{\alpha} , \bar{\beta}) \quad )
\quad ) \\
&=& (\quad (\quad (\beta_h, \beta'_k) + \psi_1( \bar{\alpha} , \bar{\beta}), \quad   S_1( \bar{\alpha} , \bar{\beta})\quad ),
\\ & & 
 ( \quad (\alpha_h, \alpha'_k) + \psi_2( \bar{\alpha} , \bar{\beta}),  \quad  S_2( \bar{\alpha} , \bar{\beta})\quad ) \quad )
\end{eqnarray*} 
where 
\begin{eqnarray*} 
\psi_1(  \bar{\alpha} , \bar{\beta}) &=& (\alpha_{h-1}-\beta_{h-1}, 0 )\quad \in \quad A= 
  \Z_q \times \Z_q \\
\psi_2(  \bar{\alpha} , \bar{\beta}) &=& (0, \beta'_{k-1}-\alpha'_{k-1} ) \quad \in \quad A= 
\Z_q \times \Z_q
\end{eqnarray*}
so that $V$ has the desired description, and explicit formulas 
for $\psi_1$ and $\psi_2$ are also obtained. 
\qed

\section{The Yang-Baxter cocycle knot invariants} \label{invsec}

\subsection{A brief review of biquandles and virtual knots} 
\label{prelimsubsec}

A {\it quandle}, $X$, is a set with a binary operation $(a, b) \mapsto a * b$
such that

(I) For any $a \in X$,
$a* a =a$.

(II) For any $a,b \in X$, there is a unique $c \in X$ such that 
$a= c*b$.

(III) 
For any $a,b,c \in X$, we have
$ (a*b)*c=(a*c)*(b*c). $

A {\it rack} is a set with a binary operation that satisfies 
(II) and (III).
Racks and quandles have been studied in, for example, 
\cite{Brieskorn,FR,Joyce,K&P,Matveev}.
A generalization of racks and quandles have been studied in several papers. 
Here we follow descriptions in \cite{BF,FJKW}.

\begin{definition} {\bf  \cite{BF,FJKW}} {\rm 
A set $X$ is called a birack if  
there is a mapping
$R: X \times X \rightarrow X \times X$ 
with the following properties.

\begin{enumerate}
\item
The map $R$ is invertible.
The inverse of $R$ is denoted by $\bar{R}: X \times X \rightarrow X\times X$. 

The images of the map $R$ are written  by
$$R(A_1, A_2)=(R_1 (A_1, A_2), R_2 (A_1, A_2) )= (A_3, A_4) , $$
where $A_i \in X$ for $i=1,2,3,4$.

\item 
For any $A_1, A_3 \in X$ there is a unique $A_2 \in X$ such that
$R_1 (A_1, A_2) = A_3$.
We say that $R_1$ is left-invertible. 

\item 
For any $A_2, A_4 \in X$ there is a unique
 $A_1 \in X$ such that $R_2 (A_1, A_2)= A_4$. 
We say that  $R_2$ is right-invertible.

\item 
$R$ satisfies the set-theoretic Yang-Baxter equation:
$$ (R \times 1) (1 \times R) (R \times 1)
=  (1 \times R)(R \times 1) (1 \times R) ,$$
where $1$ denotes the identity mapping.
\end{enumerate}

To specify the map, we also say that $(X, R)$ is a birack.

} \end{definition}

\begin{figure}[h]
\begin{center}
\mbox{
\epsfxsize=2in
\epsfbox{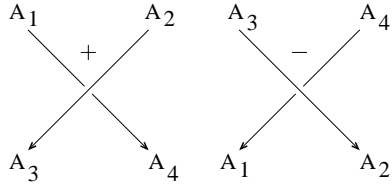} 
}
\end{center}
\caption{A coloring by birack elements }
\label{Rmatrix} 
\end{figure}

\begin{definition}{\bf \cite{BF,FJKW} } \label{biqdef} 
{\rm 
A biquandle $(X, R)$ is a birack with the following property,
called the {\em type I condition}. 

Given an element $a$ in $X$, there exists a unique $x \in X$ 
such that $x=R_1(x,a)$ and $a=R_2(x,a)$. 

} \end{definition}

\begin{remark} \label{biqrem} {\rm 
The definition of the type I condition in \cite{BF,FJKW} 
is differently formulated but they prove that their definition
is equivalent to the above condition. 
They also prove that the condition is equivalent to 
that given  $a$ in $X$, there exists a unique $y \in X$ 
such that $y=R_2(a,y)$ and $a=R_1(a,y)$.
In particular, one of these implies both.
} \end{remark} 

\begin{example} \label{biqex} {\rm
Those Yang-Baxter sets described in Example~\ref{YBEexamples} 
define biquandles.  
} \end{example}

\bigskip

\begin{figure}[h]
\begin{center}
\mbox{
\epsfxsize=3.5in
\epsfbox{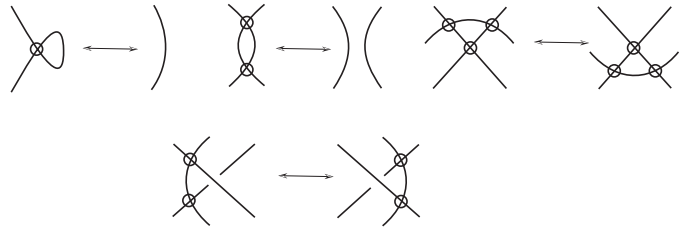} 
}
\end{center}
\caption{ Virtual Reidemeister moves  }
\label{VRMove} 
\end{figure}

A {\em virtual link  diagram\/} \cite{KauA,KauB} consists of 
generic closed curves in
${\R}^2$ such that each crossing is either a classical crossing with
over- and under-arcs, or a virtual crossing without over or under
information.
Let ${\cal VL}$ be the set of virtual link diagrams.
{\em Virtual   Reidemeister equivalence} is an equivalence relation
on ${\cal VL}$ generated by the Reidemeister moves
depicted  in Fig.~\ref{VRMove}, and ordinary Reidemeister moves.
Put $VL = {\cal VL}/_{\sim}$, where $\sim$ is the
virtual Reidemeister equivalence.  Each element of
$VL$ is called a {\em virtual link\/}.

If the given set of curves of a diagram
is connected (i.e, the diagram consists of
a single component curve), then it is called a {\em virtual
knot 
diagram\/}.
The set of virtual knot diagrams  is 
denoted by ${\cal VK}$, and
the set of equivalence classes is 
denoted by
 $VK = {\cal VK}/_{\sim}$, whose elements are called
 {\em virtual knots\/}.

\bigskip

Let $K$ be 
 an  oriented virtual link diagram.
A diagram is  immersed circles with crossing information
(over and under arcs) specified at each crossing.
The under-arcs are broken into two pieces at a crossing 
to specify the crossing information. In this case, the rest consists
of arcs, called {\it over-arcs} of the diagram.
The 
family of 
underlying immersed circles 
is 
called a {\it projection}.
When crossing points are deleted from the projection, we obtain a set of 
{\it open arcs}, denoted by ${\cal A}$. In this case,
over-arcs are also broken at a crossing. 

\begin{definition}{\rm
Let $(X, R)$ be a birack. 
Let $K$ be an oriented link diagram, ${\cal A}$ the set of its open arcs.
A map ${\cal C}: {\cal A} \rightarrow X$ is called 
a {\it coloring} of $K$ by $X$ if it satisfies the following conditions.

Let $A_i$, $i=1,2,3,4$,  be the images of open arcs near a positive 
crossing as depicted in the left of 
Fig.~\ref{Rmatrix} under the map ${\cal C}$,
so that they are elements of $X$ (regarded as being assigned to each arc).
Then it is required that 
$R(A_1, A_2)=(A_3, A_4)$. 
At a negative crossing, the elements $A_i$ as depicted in the right 
of Fig.~\ref{Rmatrix} are required to satisfy 
$\bar{R}(A_1, A_2)=(A_3, A_4)$. 
At a virtual crossing the colorings do not change for either arcs involved.

The images $A_i$ are called {\it colors}. 

} \end{definition}

\bigskip

It is proved in \cite{FJKW} that the number of colorings 
of a virtual knot diagram by a biquandle does not depend on the choice of 
a diagram, so that the number of colorings is an invariant
of virtual knots. 
In fact, it was shown that there is a one-to-one correspondence 
between the set of colorings before and after each 
Reidemeister or  virtual Reidemeister move, such that the colors
outside of a small neighborhood where the move occurs
are fixed before and after the move.

\subsection{Yang-Baxter cocycle invariants} \label{cocyinvsubsec}

In this section we define  knot invariants
by state-sum, using Yang-Baxter cocycles.
They 
 generalize
 the quandle cocycle invariant \cite{CJKLS}, 
and include it as a special case. 
The latter has been generalized to 
knotted surfaces in $4$-space, and used for topological applications 
\cite{CJKLS,SS}, so development into such directions 
are expected for the Yang-Baxter cocycle invariants as well.

 In the case of quandle cocycle invariants, we needed to 
define ``quandle condition'' \cite{CJKLS} for rack cocycles 
so that the state-sum is invariant under type I Reidemeister move. 
We need the following condition for Yang-Baxter cocycles
for this purpose. 

\begin{definition} {\rm 
Let $(X, A)$ be a biquandle, 
and $A$ be an abelian group. 
Recall in Definition~\ref{biqdef} and Remark~\ref{biqrem} 
that for any $a \in X$ there is a unique $x$ and $y \in X$ such that 
$x=R_1(x,a)$, $a=R_2(x,a)$, $y=R_2(a,y)$, and $a=R_1(a,y)$. 

A Yang-Baxter $2$-cocycle $\psi \in  Z^2_{\rm YB}(X;A)$
is said to satisfy the {\it type I condition} if 
$\psi( x,a)=0$ and $\psi(a,y)=0$ for any $a \in X$.
} \end{definition}

\bigskip

\begin{figure}[h]
\begin{center}
\mbox{
\epsfxsize=3in
\epsfbox{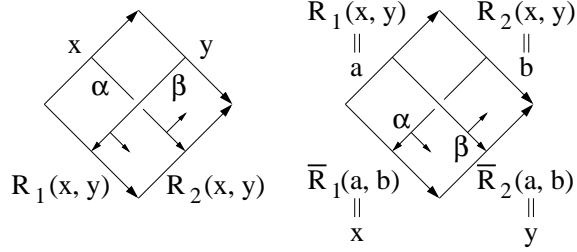} 
}
\end{center}
\caption{Crossings of classical knot diagrams }
\label{crossings} 
\end{figure}

Let $K$ be a classical knot or link diagram. Let a finite 
biquandle $(X,R)$, and a $2$-cocycle 
$\psi \in Z^2_{\rm YB}(X;A)$ 
be given,
where $A$ is an abelian group.
Let  ${\cal  C}$ 
denote a coloring ${\cal  C}: {\cal E} \rightarrow X$,
where $ {\cal E} $ denotes the set of open arcs of $K$.

Fix normals to the arcs in such a way that the (tangent, normal) 
matches the orientation of the plane, see Fig.~\ref{crossings}.
Note that the edges of the squares dual to the crossings have 
the parallel orientations to the normals, while tangents may not.

Let $\alpha$ be 
the under-arc away from which the normal to the over-arc points.
Let $\beta$ be 
the over-arc towards which the normal to the 
under-arc points. 
Let ${\cal C}(\alpha)=x$ and  ${\cal C}(\beta)=y$, see Fig.~\ref{crossings}. 
A {\it (Boltzmann) weight}, $B(\tau, {\cal C})$ 
(that depends on $\psi$),
at a  crossing $\tau$ is defined 
by $B(\tau, {\cal C}) = \psi(x,y)^{\epsilon (\tau)}$, 
where $\epsilon (\tau)=1$ or $-1$ if $\tau$ is positive or negative, 
respectively.

The {\it (Yang-Baxter) cocycle knot invariant} is defined by 
the state-sum expression  
$$
\Phi_{\rm YB} (K) = \sum_{{\cal C}}  \prod_{\tau}  B( \tau, {\cal C}). 
$$
The product is taken over all crossings of the given diagram $K$,
and the sum is taken over all possible colorings.
The values of the state-sum
are  taken to be in  the group ring $\Z[A]$ where $A$ is the coefficient 
group  written multiplicatively. 
The state-sum depends on the choice of $2$-cocycle $\psi$.

By checking all the Reidemeister moves as in \cite{CJKLS}, we obtain 

\begin{proposition} 
If $\psi \in  Z^2_{\rm YB}(X;A)$ satisfies the type I condition, 
then $\Phi_{\rm YB} (K)$ does not depend on the choice of a diagram.
Hence  $\Phi_{\rm YB} (K)$ defines a knot invariant. 
\end{proposition}

\begin{proposition} \label{obstprop}
Let $\psi \in   Z^2_{\rm YB}(X;N)$ be an obstruction
 cocycle as in Equality~(\ref{defining2}).
Then the value of the  Yang-Baxter cocycle knot invariant 
$\Phi_{\rm YB}(K)$ with $\psi$ is a positive integer for any 
classical or virtual knot or a link $K$.
\end{proposition}
{\it Proof.\/} 
The proof is similar to the one found in \cite{CJKLS}.
Choose and fix any coloring of a given knot diagram by a biquandle.
Suppose $f \in  Z^2_{\rm YB}(X;A)$ be as in Equality (\ref{defining2}). 
Since $i: N \rightarrow G$ is injective, compute the state-sum in 
$\Z [G] $ which contains $\Z [N] $. 
Then the weight  assigned at a crossing, say, $\tau$
 is $i \psi(x, y)^{\epsilon} $,
where $\epsilon=1$ or $-1$ depending on whether $\tau$ 
is positive or negative, respectively.
Assign the terms $sf(x)^{\epsilon} $, $sf(y)^{\epsilon} $,
 $sf(R_1(x,y))^{-{\epsilon} }$, and $sf(R_2(x,y))^{-{\epsilon} }$
on the strings near the crossing  $\tau$ colored with 
$x$, $y$, $R_1(x,y)$, and $R_2(x,y)$, respectively.
Then the state-sum is computed by using these weights.
However, the weights assigned to the two ends of each arc cancel, 
giving $1$ as the contribution to the state-sum from this coloring.
\qed

\subsection{Computations} \label{compsubsec}

In \cite{BF}, it was pointed out that if a matrix
$\displaystyle R= \left[ \begin{array}{cc} A & B \\ C & D \end{array} \right] $
satisfies the YBE, then so does 
$\displaystyle R= \left[ \begin{array}{cc} A & u B \\ u^{-1} C & D \end{array} \right] $
for any invertible element $u$ 
in the center of the ring used.
In particular, consider such variants 
$\displaystyle R= \left[ \begin{array}{cc} 1-s & u s \\ u^{-1} t & 1-t \end{array} \right] $ of those in Example~\ref{YBEexamples}.
This  matrix defines 
a biquandle structure on $ \Z_q$ for a 
positive integer $q$, 
if $(1-s)(1-t)\equiv 0 \pmod{q}$ 
and $s,t,u$ are 
all invertible. 
In this section we present some calculations for this type of biquandles.

\begin{example} {\rm 
We computed the number of colorings 
of Kishino-type knots (see, for example, \cite{BF,FJKW})
 for $q=15$, 
 $s=4$, $t=11$, and  $u=2,4,7,8,11,13,14$.
The Kishino's knot is a connected sum of diagrams of the unknot,
and we took three 
variants usually considered 
($K_1$ through $K_3$)
and further variants 
($K_4$ through $K_6$), 
represented by 
the closed braid form of the following braid words.
We use the usual convention $\sigma_i$ for a standard 
braid generator between $i$th and $(i+1)$st strings, 
and the notation $v_i$ for virtual crossing between them.
The results are summarized in the table below.
We included these examples for noting that 
the numbers of colorings alone distinguish some of these from others,
as well as from the unknot.

\begin{eqnarray*}
K_1 & = & \sigma_1 v_1 \sigma_1^{-1} \sigma_2 \sigma_1 v_1
 \sigma_1^{-1} \sigma_2^{-1} \\
K_2 & = & \sigma_1 v_1 \sigma_1^{-1} \sigma_2 \sigma_1^{-1} v_1
 \sigma_1 \sigma_2^{-1} \\
K_3 & = & \sigma_1^{-1} v_1 \sigma_1 \sigma_2 \sigma_1^{-1} v_1
 \sigma_1 \sigma_2^{-1} \\
K_4 & = &\sigma_1 \sigma_1 v_1  \sigma_1^{-1} \sigma_1^{-1}
\sigma_2 \sigma_1 \sigma_1 v_1  \sigma_1^{-1}  \sigma_1^{-1} \sigma_2^{-1} \\
K_5 & = &\sigma_1 \sigma_1 v_1  \sigma_1^{-1} \sigma_1^{-1}
\sigma_2 \sigma_1^{-1} \sigma_1^{-1} v_1  \sigma_1 \sigma_1 \sigma_2^{-1} \\
K_6 & = &\sigma_1^{-1} \sigma_1^{-1} v_1  \sigma_1 \sigma_1
\sigma_2 \sigma_1^{-1} \sigma_1^{-1} v_1  \sigma_1  \sigma_1 \sigma_2^{-1} \\
\end{eqnarray*}

\begin{table}[h] 
\begin{center}
{\begin{tabular}{||l||r|r|r|r|r|r|r|r||}\hline \hline
 $u$  & $K_1$ & $K_2 $ &   $K_3 $ & $K_4 $ & $K_5 $ & $K_6 $ 
\\ \hline \hline 
2&225&15&
75&15&15&
45
\\ \hline 
4&15&15&
45&45&15&
75
\\ \hline 
7&75&15&
225&45&15&
15
\\ \hline
8&45&15&
15&75&15&
225
\\ \hline
11&45&15&
15&75&15&
45
\\ \hline
13&15&15&
45&225&15&
75
\\ \hline
14&45&15&
15&15&15&
225
\\ \hline \hline
\end{tabular} } \end{center}
\caption{A table of numbers of colorings of Kishino-type knots}
\label{kishinotable}
\end{table}

} \end{example}

\begin{example}{\rm 
Let $X$ be the biquandle determined by the matrix 
$\displaystyle R= \left[ \begin{array}{cc} 0 & 2 \\ 1 & 2 \end{array} \right] $ in $\Z_3$ as in the preceding example,
where $q=3$, 
$s=1$, $t=2$, and $u=2$.  
Then we obtain a non-trivial $2$-cocycle 
$f \in Z_{\rm YB}^2(X; \Z_3)$ with the following values:
$$ f(1,0)=q_1, \ f(2,2)=q_2, \ f(1,1)=q_3, f(2,0)=-q_1, \\
f(0,2)=q_1 - q_3, \ f(0,1)=-q_1 - q_2, $$
where $q_i$, $i=1,2,3$, are arbitrary choices of numbers modulo $3$.
The values for unspecified evaluations are zeros.
This cocycle was obtained by {\it Maple} and confirmed independently 
by {\it Mathematica}.  
This small biquandle of three elements already takes  non-trivial 
values for the cocycle invariant on virtual knots. We obtain, for example, 
the following by direct calculations:
\begin{equation*}
\Phi _{\rm YB}({\rm cl}(\sigma_1^n v_1))=\begin{cases}
3 \hspace{1in} \text{if} & n \equiv 0\; (\rm mod\; 3)\\
1 + \xi^{q_1} + \xi^{-q_1} & n \equiv 1, 2\; (\rm mod \;3)  .
\end{cases}
\end{equation*}
 Here, $\xi$ denotes a multiplicative generator of the coefficient group $\Z_3$.
}\end{example}

\begin{example}{\rm 
Let $X$ be the biquandle determined by the matrix 
$\displaystyle R= \left[ \begin{array}{cc} 0 & -1 \\ 1 & 2 \end{array} \right] $ in $\Z_4$ as 
before,  where $q=4$, 
$s=1$, $t=-1$, and $u=-1$.  
Then we obtain a non-trivial $2$-cocycle 
$f \in Z_{\rm YB}^2(X; \Z_4)$ with the following values:
$$ f(0,1)=f(1,1)=f(1,2)=f(3,3)=1,  \ f(0,2)=2, \ 
f(1,0)=f(2,1)=f(3,0)=f(3,2)=3. $$
The values for unspecified evaluations are zeros.
This cocycle was again obtained by {\it Maple} and confirmed independently 
by {\it Mathematica}.

Let $\xi$ denote a generator of the coefficient group $\Z_4$.
For the torus links (the closure ${\rm cl}(\sigma_1^n)$
of  $2$-braids,  $n \in \Z$), 
we obtain the following formula.
\begin{equation*}
\Phi _{\rm YB}({\rm cl}(\sigma_1^n))=\begin{cases}
4 \hspace{1in} \text{if} & n \equiv 1\; (\rm mod\; 2)\\
4+4\xi^2 & n \equiv 2\; (\rm mod \;4)\\
8+8\xi^3 & n \equiv 4\; (\rm mod \; 16)\\
8+8\xi^2 & n \equiv 8\; (\rm mod \; 16)\\
8+8\xi & n \equiv 12\; (\rm mod \;16)\\
16 & n \equiv 0\; (\rm mod \; 16)  .
\end{cases}
\end{equation*}
{\it Proof.\/} 
Since $R^4$ is the identity matrix, and the coefficient group is $\Z_4$, it is clear that
the cocycle invariant has period $16$:
$\Phi _{\rm YB}({\rm cl}(\sigma_1^{n+16}))=\Phi _{\rm YB}({\rm cl}(\sigma_1^{n}))$
for any $n \in \Z$.
When $n \equiv 1\; (\rm mod\; 2)$, the closed braid ${\rm cl}(\sigma_1^{n})$ 
is colored if and only if the colors at the top strings are of the form 
$(a, -a)$, $a \in X$. This can be seen from the matrix $R^{n}$. 
For each such coloring the contribution at each crossing is zero,  so that  $\Phi _{\rm YB}({\rm cl}(\sigma_1^n))=4$.
Now if $n \equiv 2\; (\rm mod\; 4)$ then the closed braid is colored if and only if the colors at the top arcs  are  of the form $(a,-a)$ or $(a,2-a)$, $a \in X$.  
Any coloring of the form $(a,2-a)$ gives a total contribution of two.  
This implies that $\Phi_{\rm YB} ({\rm cl}(\sigma_1^n))=4+4\xi^2$.  
There are 16 possible colorings if $ n \equiv 4\; (\rm mod \; 16)$, for all pairs $(a,b)$
of top arcs.  
For the colors of the form $(a, -a)$, the contribution is zero. For $(a, 2-a)$,
the contribution is doubled compared to the case $n \equiv 2\; (\rm mod\; 4)$,
so that the contribution is also zero.
The remaining colors contribute $3$, so that 
$\Phi_{\rm YB} ({\rm cl}(\sigma_1^n))=8+8\xi^3$.
The contributions of these remaining colors get doubled, tripled, and quadrupled,
respectively, for the cases $n \equiv 8, 12,$ and $0 \mod 16$.   
\qed 

In particular,
$\Phi_{\rm YB}({\rm cl}(\sigma_1^4))  =  8 + 8 \xi^3 $ and 
$\Phi_{\rm YB}({\rm cl}(\sigma_1^{-4}) )  =  8 + 8 \xi$ so that 
 the invariant detects the chirarity.
 
The invariant can take other  values for some 
virtual knots, such as:
\begin{equation*}
\Phi _{\rm YB}({\rm cl}(\sigma_1 v_1)^n)=\begin{cases}
3 + \xi + \xi^2 + 3\xi^3 \hspace{1in} \text{if} & n \equiv 1\; (\rm mod\; 4)\\
6+6\xi^2+2\xi+2\xi^3 & n \equiv 2\; (\rm mod \;4)\\
3 + 3\xi + \xi^2 + \xi^3  & n \equiv 3\; (\rm mod\; 4)\\
12+4\xi^2 & n \equiv 4\; (\rm mod \; 8)\\
16 & n \equiv 0\; (\rm mod \; 8)  .
\end{cases}
\end{equation*}

The invariant can be non-trivial for links with zero linking number.
The Borromean rings, a three-component link with pairwise linking number 
zero and represented as the closed braid 
of $(\sigma_1^{-1} \sigma_2 )^3$, has the value 
$16 + 48\xi^2$.
} \end{example}

These calculations suggest that these invariants are quite non-trivial
and deserve further study.


\begin{thebibliography}{99}
\setlength{\itemsep}{-3pt}

\bibitem{BF} Bartholomew, A.; Fenn, R.,
{\it Quaternionic invariants of virtual knots and links}, preprint
found at http://www.maths.sussex.ac.uk///Staff/RAF/Maths/.


\bibitem{Brieskorn} Brieskorn, E., 
{\it Automorphic sets and singularities,}
Contemporary math., 78 (1988), 45--115.

\bibitem{CENS} 
 Carter, J.S.; Elhamdadi, M.; Nikiforou, M.A.; Saito, M.,
{\it Extensions of quandles and cocycle knot invariants,}
  J. Knot Theory Ramifications,  12  (2003),  725--738. 
  
\bibitem{CJKLS}
 Carter, J.S.; Jelsovsky, D.; Kamada, S.; Langford, L.; Saito, M.,
{\it Quandle cohomology and state-sum invariants
of knotted curves and surfaces,}
  Trans. Amer. Math. Soc.  355  (2003),   3947--3989.
  
\bibitem{CJKS1}
 Carter, J.S.; Jelsovsky, D.; Kamada, S.; Saito, M.,
{\it Computations of quandle cocycle invariants of
 knotted curves and surfaces,}
 Advances in math., 157 (2001) 36-94. 

\bibitem{CJKS2}
 Carter, J.S.; Jelsovsky, D.; Kamada, S.; Saito, M.,
{\it Quandle homology groups, their betti numbers, and virtual knots,}
 J. of Pure and Applied Algebra,
  157,  (2001), 135-155.

\bibitem{SSS2}
 Carter, J.S.; Kamada, S.; Saito, M.,
{\it Geometric interpretations of quandle homology,}
 J. of Knot Theory and its Ramifications, 10, no. 3 (2001), 345-386.


\bibitem{CS:book} Carter, J.S.; Saito, M.,
{\it Knotted surfaces and their diagrams,}
the American Mathematical Society,  1998.

\bibitem{Diamond}
 Carter, J.S.; Saito, M.,
{\it Diagrammatic invariants of knotted curves and surfaces,}
unpublished note written in 1992
with limited circulation,  on which this paper is based. 



\bibitem{ESG} Etingof, P.; Soloviev, A.; Guralnick, R.,
{\it  Indecomposable set-theoretical
solutions to the quantum Yang-Baxter equation on a set with a prime number of elements},  J. Algebra 242
(2001), no. 2, 709--719.

\bibitem{ESS} Etingof, P.; Schedler, T.; Soloviev, A.,
{\it Set-theoretical solutions of the quantum Yang-Baxter equation},
 Duke Math. J. 100 (1999), no. 2, 169--209. 

\bibitem{FJKW} 
Fenn, R.; Jordan-Santana, M.;  Kauffman, L.; Wraithe, G.,
{\it Biracks and virtual links}, 
preprint found at http://www.maths.sussex.ac.uk///Staff/RAF/Maths/.

\bibitem{FR}   Fenn, R.; Rourke,  C.,
\textit{Racks and links in codimension two.}
Journal of Knot Theory and Its Ramifications Vol. 1 No. 4 (1992), 343-406.

\bibitem{FRS0}
Fenn, R.; Rourke, C.; Sanderson, B.,
{\it An introduction to species and the rack space}, 
Topics in knot theory (Erzurum, 1992), 33--55,
NATO Adv. Sci. Inst. Ser. C Math. Phys. Sci., 399,
Kluwer Acad. Publ., Dordrecht, 1993.

 \bibitem{FRS1} 
Fenn, R.; Rourke, C.; Sanderson, B., {\it Trunks and classifying spaces,}
Appl. Categ. Structures 3 (1995), no. 4, 321--356.


\bibitem{FRS2} Fenn, R.; Rourke, C.; Sanderson, B.,
{\it James bundles and applications,} preprint found at
http://www.maths.warwick.ac.uk/${ }^{\sim}$bjs/ .

\bibitem{FoxTrip} Fox, R.H., 
{\it A quick trip through knot theory,}
in Topology of $3$-Manifolds, 
Ed. M.K. Fort Jr., Prentice-Hall (1962) 120--167.

\bibitem{GPV} Goussarov, M.; Polyak, M.; Viro, O.,
{\it  Finite-type invariants of classical and
virtual knots},  Topology 39 (2000), no. 5, 1045--1068.

\bibitem{Jones} 
Jones, V.F.R.,
{\it Hecke algebra representations of braid groups and link polynomials},
{Ann. of Math.} {126} (1989),   335--388.



\bibitem{Joyce} Joyce, D.,
{\it A classifying invariant of knots, the knot quandle,}
J. Pure Appl. Alg., 23, 37--65.

\bibitem{KK}  Kamada, N.; Kamada, S.,
{\it  Abstract link diagrams and virtual knots}, J. Knot
Theory Ramifications 9 (2000),  93--106. 


\bibitem{K&P} L. H. Kauffman, {\it Knots and Physics},
World Scientific, Series on knots and everything, vol. 1, 1991.


\bibitem{KauA} L. H. Kauffman,
{\it Virtual knots},
talks at MSRI Meeting in January 1997 and AMS Meeting at
University of Maryland, College Park in March 1997.

\bibitem{KauB} L. H. Kauffman,
{\it Virtual Knot Theory},
European J. Combin. {\bf 20} (1999), 663--690.

\bibitem{KR} Kauffman, L.H., and Radford, D.E.,
  {\it  Bi-oriented quantum algebras, and a generalized
Alexander polynomial for virtual links}, 
Diagrammatic morphisms and applications (San Francisco, CA, 2000),   113--140, Contemp. Math., 318, Amer. Math. Soc., Providence, RI,  2003. 

 \bibitem{Kim} Kim, Se-Goo, 
{\it  Virtual knot groups and their peripheral structure},  J. Knot Theory
Ramifications 9 (2000), 797--812.

  \bibitem{LYZ}    Lu, J-H.; Yan, M.; Zhu, Y-C,
{\it  On the set-theoretical Yang-Baxter
equation,}  Duke Math. J. 104 (2000), no. 1, 1--18. 

\bibitem{Matveev} Matveev, S., 
{\it Distributive groupoids in knot theory,} (Russian) Mat. Sb. (N.S.)
119(161) (1982), no. 1, 78--88, 160.
             

\bibitem{RS} Rourke, C., and Sanderson, B.,
{\it There are two $2$-twist-spun trefoils,}
preprint at 
 http://xxx.lanl.gov/abs/math.GT/0006062 .


\bibitem{SS} Satoh, S.; Shima, A.,
{\it The $2$-twist spun trefoil has the triple point number four,}
Preprint, to appear in Trans. Amer. Math. Soc.

\bibitem{Saw} Sawollek, J.,
{\it On Alexander-Conway polynomials for virtual knots and links,}
preprint, arXiv:math.GT/9912173.


\bibitem{SilWil} Silver, D., and  Williams, S.,
{\it   Alexander groups and virtual links},  J. Knot Theory
Ramifications 10 (2001),  151--160. 

\bibitem{DanSusan} Silver, D., and  Williams, S.,
{\it A generalized Burau representation for string links, } 
Pacific J. Math. 197 (2001),  241--255.

 \bibitem{Solo}   Soloviev, A.,
{\it  Non-unitary set-theoretical solutions to the quantum Yang-Baxter
equation,} Math. Res. Lett. 7 (2000), no. 5-6, 577--596.

\bibitem{Turaev} Turaev, V.,
{\it The Yang-Baxter equation and invariants of links,}
Invent. math. 92 (1988) 527--553.

\bibitem{Rose} Roseman, D.,
{\it Reidemeister-type moves for surfaces in four dimensional space, } 
in Banach Center Publications 42 (1998) Knot theory, 347--380.


\bibitem{Rolf} Rolfsen, D.,
{\it Knots and Links.}  Publish or Perish Press, (Berkley 1976).

\bibitem{Wada} Wada, M.,
{\it Group invariants of links,}
Topology 31 (1992), 399-406. 


\end{thebibliography}
\end{document}